\documentclass[12pt]{article}%
\usepackage{amsmath}
\usepackage{amsfonts}
\usepackage{amssymb}
\usepackage{graphicx}
\providecommand{\U}[1]{\protect\rule{.1in}{.1in}}

\begin{document}

\title{5-Engel Lie algebras II}
\author{Michael Vaughan-Lee}
\date{October 2024}
\maketitle

\begin{abstract}
In my article $5$\emph{-Engel Lie algebras} published on the arXiv in 2023 I
proved that $5$-Engel Lie algebras of characteristic zero or prime
characteristic $p>7$ are nilpotent of class at most $11$. In this note I
investigate the ideal Id$(x)$ generated by an element $x$ in a $5$-Engel Lie
algebra. Traustason [5] has given an example of a $3$-Engel Lie algebra of
characteristic 2 in which Id$(x)$ does not have to be nilpotent, and in this
note I give an example of a $5$-Engel Lie algebra of characteristic $3$ in
which Id$(x)$ does not have to be nilpotent. For all primes $p>3$, I show that
Id$(x)$ is nilpotent in $5$-Engel Lie algebras of characteristic $p$, and I
obtain explicit (best possible) bounds on the nilpotency class.

\end{abstract}

\section{Introduction}

A Lie algebra $L$ is said to satisfy the $n$-Engel identity if%
\[
\lbrack x,\underset{n}{\underbrace{y,y,\ldots,y}}]=0
\]
for all $x,y\in L$. It follows from Zel'manov's solution of the restricted
Burnside problem \cite{zelmanov91}, \cite{zelmanov91b} that $n$-Engel Lie
algebras are locally nilpotent (for all $n$). And in a very deep theorem
Zel'manov \cite{zelmanov87} also showed that $n$-Engel Lie algebras over a
field of characteristic zero are globally nilpotent (for all $n$). Note that
if $n$-Engel Lie algebras over fields of characteristic zero are nilpotent of
class at most $c$ then $n$-Engel Lie algebras of prime characteristic $p$ will
also be nilpotent of class at most $c$ for all sufficiently large $p$. In
\cite{vlee5eng} I proved that 5-Engel Lie algebras of characteristic zero, or
of prime characteristic $p>7$, are nilpotent of class at most 11. On the other
hand Razmyslov \cite{Razmyslov71} has shown that there are non-soluble (and
hence non-nilpotent) $(p-2)$-Engel Lie algebras of characteristic $p$ for all
primes $p\geq5$. For $p=5$ this result was originally proved by Bachmuth and
Mochizuki \cite{bachmuth70}. So Razmyslov's construction gives non-nilpotent
5-Engel Lie algebras of characteristic 5 and 7. Also, it is easy to construct
non-nilpotent 3-Engel Lie algebras of characteristic 2 and non-nilpotent
4-Engel Lie algebras of characteristic 3 (see \cite{vlee5eng}).

In this note I investigate Id$(x)$ for elements $x$ in 5-Engel Lie algebras.
Traustason \cite{Traustason93} gives an example of a 3-Engel Lie algebra of
characteristic 2 in which the ideal generated by an element does not have to
be nilpotent, and clearly his example is also 5-Engel. In Section 2, I give an
example of a 5-Engel Lie algebra of characteristic 3 in which the ideal
generated by an element does not have to be nilpotent. And in Section 3$\,$ I
describe computer calculations which show that in characteristic 5 the ideal
generated by an element in a 5-Engel Lie algebra is nilpotent of class at most
7; in characteristic 7 and 19 it is nilpotent of class at most 6; and in
characteristic $p>7$ ($p\neq19$) it is nilpotent of class at most 5. These
results are best possible.

\section{An example in characteristic 3}

Let $A$ be the class two Lie algebra over GF(3) generated by $a_{1}%
,a_{2},a_{3},\ldots$ subject to the relations $[a_{i},a_{j}]=0$ for all
$i,j>1$. So $[A,A]$ is central in $A$ and is infinite dimensional with basis
$[a_{2},a_{1}],[a_{3},a_{1}],\ldots$. Let $B$ be a one dimensional Lie
algebra over GF(3) generated by the element $b$, and let $L$ be the Wreath
product of $B$ with $A$. Let $C$ be the ideal of $L$ generated by $b$. So
$L=A+C$, and $C$ is an abelian ideal with basis consisting of all distinct
elements of the form%
\begin{equation}
\lbrack b,a_{i},a_{j},\ldots,a_{k},[a_{r},a_{1}],[a_{s},a_{1}],\ldots ,\lbrack
a_{t},a_{1}]]\tag{*}%
\end{equation}
with $i\leq j\leq\ldots\leq k$, $r\leq s\leq\ldots\leq t$. (The sequences
$i,j,\ldots,k$ and $r,s,\ldots,t$ can be empty.)

Let $I$ be the ideal of $L$ spanned by all elements of the form%
\[
\lbrack b,a_{1},a_{1},a_{1},a_{i},a_{j},\ldots,a_{k},[a_{r},a_{1}%
],[a_{s},a_{1}],\ldots ,\lbrack a_{t},a_{1}]].
\]
(So $I$ is the ideal of $L$ generated by $[b,a_{1},a_{1},a_{1}]$.) And let $J
$ be ideal of $L$ spanned by all elements of the form (*) with the property
that there exists $n>1$ such that two (or more) elements of the sequence%
\[
i,j,\ldots,k,r,s,\ldots,t
\]
are equal to $n$. We let $M=L/(I+J)$. The ideal $C/(I+J)$ in $M$ has a basis
consisting of all elements of the form (*) which have degree at most 1 in any
$a_{i}$ with $i>1$, and where the sequence $i,j,\ldots,k$ has at most two
entries equal to 1. Clearly the ideal of $M$ generated by $a_{1}+I+J$ is not
nilpotent. We show that $M$ is 5-Engel.

Since $[M,M,M]\leq C/(I+J)$, to show that $M$ is 5-Engel it is sufficient to
show that $[c,a,a,a]\in I+J$ for all elements $c\in C$ and for all $a\in A$.
And to this end it is sufficient to show that if $x,y,z$ are three distinct
elements of the basis $a_{1},a_{2},\ldots,[a_{2},a_{1}],[a_{3},a_{1}],\ldots$
for $A$, and if $c\in C$ then

\begin{enumerate}
\item $[c,x,x,x]\in I+J$,

\item $[c,x,x,y]+[c,x,y,x]+[c,y,x,x]\in I+J$,

\item $[c,x,y,z]+[c,x,z,y]+[c,y,x,z]+[c,y,z,x]+[c,z,x,y]+[c,z,y,x]\in I+J$.
\end{enumerate}

First consider the case (1). Clearly if $x$ is a basis element other than
$a_{1}$ then $[c,x,x,x]\in J$ for all $c\in C$. So consider the case when
$x=a_{1}$ and $c$ is a basis element for $C$ of the form (*).%
\begin{align*}
& [[b,a_{i},a_{j},\ldots,a_{k},[a_{r},a_{1}],[a_{s},a_{1}],\ldots ,\lbrack
a_{t},a_{1}]],a_{1},a_{1},a_{1}]\\
& =[b,a_{i},a_{j},\ldots,a_{k},a_{1},a_{1},a_{1},[a_{r},a_{1}],[a_{s}%
,a_{1}],\ldots ,\lbrack a_{t},a_{1}]]\\
& =[b,a_{i},a_{j},\ldots,a_{1},a_{1},a_{1},a_{k},[a_{r},a_{1}],[a_{s}%
,a_{1}],\ldots ,\lbrack a_{t},a_{1}]]\\
& +3[b,a_{i},a_{j},\ldots,a_{1},a_{1},[a_{k},a_{1}],[a_{r},a_{1}],[a_{s}%
,a_{1}],\ldots ,\lbrack a_{t},a_{1}]]\\
& =[b,a_{i},a_{j},\ldots,a_{1},a_{1},a_{1},a_{k},[a_{r},a_{1}],[a_{s}%
,a_{1}],\ldots ,\lbrack a_{t},a_{1}]].
\end{align*}
Continuing in this way, collecting the three occurrences of $a_{1}$ to the
left, we eventually obtain%
\[
\lbrack b,a_{1},a_{1},a_{1},a_{i},a_{j},\ldots,a_{k},[a_{r},a_{1}%
],[a_{s},a_{1}],\ldots ,\lbrack a_{t},a_{1}]]\in I.
\]

Next consider case (2). Clearly $[c,x,x,y]+[c,x,y,x]+[c,y,x,x]\in J$ unless
$x=a_{1}$. And collecting $a_{1}$ to the left in the case when $x=a_{1}$ we
see that%
\begin{align*}
& \lbrack c,x,x,y]+[c,x,y,x]+[c,y,x,x]\\
& =3[c,a_{1},a_{1},y]+3[c,a_{1},[y,a_{1}]]\\
& =0.
\end{align*}

Finally consider case (3). If none of $x,y,z$ are equal to $a_{1}$ then
$x,y,z$ all commute, and%
\begin{align*}
& [c,x,y,z]+[c,x,z,y]+[c,y,x,z]+[c,y,z,x]+[c,z,x,y]+[c,z,y,x]\\
& =6[c,x,y,z]\\
& =0.
\end{align*}
So suppose that $x=a_{1}$. Since $y,z$ are elements of the basis for $A$ which
are not equal to $a_{1}$ we see that $[y,z]=0$. Collecting $a_{1}$ to the left
as before we see that%
\begin{align*}
& [c,x,y,z]+[c,x,z,y]+[c,y,x,z]+[c,y,z,x]+[c,z,x,y]+[c,z,y,x]\\
& =[c,a_{1},y,z]+[c,a_{1},z,y]+[c,y,a_{1},z]+[c,y,z,a_{1}]+[c,z,a_{1}%
,y]+[c,z,y,a_{1}]\\
& =6[c,a_{1},y,z]+3[c,y,[z,a_{1}]]+3[c,z,[y,a_{1}]]\\
& =0.
\end{align*}

\section{5-Engel Lie algebras of characteristic $p>3$}

We show that if $L$ is a 5-Engel Lie algebra of characteristic $p>3$, and if
$x\in L$, then the ideal of $L$ generated by $x$ is nilpotent, and we obtain
bounds for the nilpotency class of Id$_{L}(x)\,$, with the bounds depending on
$p$. We make use of a beautiful idea introduced by Graham Higman
\cite{higman56a} in his solution of the restricted Burnside problem for
exponent 5. The associated Lie rings of groups of exponent 5 have
characteristic 5, and satisfy the 4-Engel identity. Higman proved that if
$x\in F$, where $F$ is a 4-Engel Lie algebra of characteristic 5, then
Id$_{F}(x)$ is nilpotent of bounded class. If the bound on the class is $c$,
then it follows that the associated Lie ring of the $d$ generator Burnside
group of exponent 5, $B(d,5)$, is nilpotent of class at most $cd$, and it
follows that finite quotients of $B(d,5)$ have class at most $cd$. (Higman did
not obtain an explicit bound for the nilpotency class of Id$_{F}(x)$, but in
fact computer calculations show that the class is at most 6.)

Higman's first step in his proof was to use a very clever hand argument to
show that if $L$ is a 4-Engel Lie algebra of characteristic 5 generated by
elements $x,a_{1},a_{2},a_{3},\ldots$ where $[a_{i},a_{j}]=0$ for all $i,j$,
then $L$ is nilpotent. (Again, he did not obtain an explicit bound on the
class of $L$, but computer calculations show that the class is 12.) 
Since $L$ is nilpotent, it follows
that Id$_{L}(x)$ is nilpotent, of class $c$ say. Higman then went on to deduce
that this result implies that if $F$ is \emph{any} 4-Engel Lie algebra of
characteristic 5, and if $x\in F$, then Id$_{F}(x)$ is nilpotent of class $c$.
This deduction relies on the fact that in characteristic 5 the 4-Engel
identity is equivalent to the multilinear identity%
\[
\sum_{\sigma\in\text{Sym}(4)}[y,x_{1\sigma},x_{2\sigma},x_{3\sigma}%
,x_{4\sigma}]=0.
\]
This implies that if $F$ is a free 4-Engel Lie algebra of characteristic 5,
freely generated by $x,a_{1},a_{2},a_{3},\ldots$ then $F$ is
\emph{multigraded}. Every Lie product of the free generators of $F$ can be
assigned a multidegree $d$, where $d$ specifies the degree of the product in
each free generator (i.e. the number of occurrences of the free generator in
the product). If for each multidegree $d$ we let $W_{d}$ be the subspace of
$F$ spanned by Lie products with multidegree $d$, then $F$ is the direct sum
of the subspaces $W_{d}$.

So suppose that Id$_{L}(x)$ is nilpotent of class $c$. We let $F$ be a free
4-Engel Lie algebra of characteristic 5, freely generated by
$x,a_{1},a_{2},a_{3},\ldots$, and we let $u$ be a Lie
product of the free generators which has total degree $c+1+k$, with
$c+1$ entries $x$ and $k$ entries from the set $\{a_{1},a_{2},\ldots\}$.
We show by induction on $k$ that $u=0$. This proves that Id$_{F}(x)$ is
nilpotent of class $c$. The case $k=0$ is trivial, and the case $k=1$ it is also
trivial since then $u$ must be a linear multiple of $[a_{i},\underset
{c+1}{\underbrace{x,x,\ldots,x}}]$ for some $i$. (Clearly we must have
$c\geq3$.) We can assume by induction that any Lie product of the free
generators of $F$ with $c+1$ entries $x$ and with total degree less than
$c+1+k$ is zero. Since Id$_{L}(x)$ is nilpotent of class $c$, $u$ must lie in
the ideal of $F$ generated by the Lie products $[a_{i},a_{j}]$, so we can
express $u$ as a linear combination of Lie products $[a_{i},a_{j},\ldots]$
where all the unspecified entries are free generators of $F$. Since $F$ is
multigraded we can assume that all the Lie products in this linear combination
have $c+1$ entries $x$ and total degree $c+1+k$. Now consider a term
$v=[a_{i},a_{j},\ldots]$ in this linear combination and let $a_{r}$ be a free
generator of $F$ which does not appear in $v$. Let $\theta$ be the
endomorphism of $F$ mapping $x$ to $x$, mapping $a_{n}$ to $a_{n}$ for $n\neq
r$, and mapping $a_{r}$ to $[a_{i},a_{j}]$. Then%
\[
v=[a_{i},a_{j},\ldots]=[a_{r},\ldots]\theta.
\]
By induction $[a_{r},\ldots]=0$, and so $v=0$. Since this is true for all $v
$, we see that $u=0$ as claimed.

Our aim is to use Higman's argument in 5-Engel Lie algebras of characteristic
$p>3$. We let $F$ be the free 5-Engel Lie algebra over a field $K$ of
characteristic $p$, freely generated by $x,a_{1},a_{2},\ldots$. And we let $L$ be
a 5-Engel Lie algebra over $K$ generated by $x,a_{1},a_{2},\ldots$ satisfying
the additional relations $[a_{i},a_{j}]=0$ for all $i,j$. For every prime
$p>3$ we obtain an explicit bound for the class of $L$, the class of
Id$_{L}(x)$, and a bound on the number of entries $a_{i}$ in a non-zero Lie
product in $L$. We then use Higman's argument to bound the class of
Id$_{F}(x)$. These bounds are given in the following table. The first four
lines in the table give the bounds for $p=5,7,19,31$, and the last line gives
the bounds for all other primes greater than 7.%
\[%
\begin{tabular}
[c]{|c|c|c|c|}\hline
$p$ & Class of $L$ & Class of Id$_{L}(x)$ & Maximum number of entries $a_{i}%
$\\\hline
5 & 13 & 7 & 7\\\hline
7 & 11 & 6 & 6\\\hline
19 & 9 & 6 & 5\\\hline
31 & 10 & 5 & 5\\\hline
& 9 & 5 & 5\\\hline
\end{tabular}
\]

For the primes 5,7,19,31 it is a straightforward calculation using the
nilpotent quotient algorithm \cite{havasnvl90} to establish these results. For
example, to establish these claims when $p=5$ we let $M$ be the 5-Engel Lie
algebra over GF$(5)$ generated by $x,a_{1},a_{2},\ldots,a_{8}$ satisfying the
relations $[a_{i},a_{j}]=0$ for all $i,j$. We also suppose that any Lie
product of the generators of $M$ with multidegree $(d_{x},d_{1},d_{2}%
,\ldots,d_{8})$ in $x,a_{1},a_{2},\ldots,a_{8}$ is zero if $d_{x}>8$ or
\thinspace$d_{i}>1$ for any $i$. Then the nilpotent quotient algorithm shows
that $M$ has class 13, and that any Lie product of the generators in $M$ is
zero if it has multidegree $(d_{x},d_{1},d_{2},\ldots,d_{8})$ with $d_{x}>7$,
or $\sum_{i=1}^{8}d_{i}>7$. This establishes our claims for $p=5$.

The bounds for the primes 7, 19, 31 are obtained in the same way.

There is a slight difficulty in applying Higman's argument that the class of
Id$_{F}(x)$ is 7 when $p=5$, since when $p=5$ the 5-Engel identity is not
equivalent to the multilinear identity%
\[
\sum_{\sigma\in\text{Sym}(5)}[y,x_{1\sigma},x_{2\sigma},x_{3\sigma}%
,x_{4\sigma},x_{5\sigma}]=0.
\]
But when $p=5$ the 5-Engel identity does follow from this multilinear identity
together with relations $[y,z,z,z,z,z]=0$ where $z$ is a Lie product of the
generators. So $F$ is still multigraded, even when $p=5$, and Higman's
argument is still valid.

To establish the last line of the table above we need a different approach,
since we need to deal with infinitely many primes simultaneously. But I was
able to use my \textsc{Magma} \cite{boscan95} implementation of the nilpotent
quotient algorithm, which can compute finite dimensional multigraded Lie
algebras over the rationals $\mathbb{Q}$. With this program I was able to
establish the following results in various subalgebras of $L$.

\begin{itemize}
\item All Lie products of multidegree $(6,1,1)$ in $\langle x,a_{1}%
,a_{2}\rangle$ are zero provided $p\neq2,3,5,7$.

\item All Lie products of multidegree $(6,1,1,1)$ in $\langle x,a_{1}%
,a_{2},a_{3}\rangle$ are zero provided $p\neq2,3,5,7,19$.

\item All Lie products of multidegree $(6,1,1,1,1)$ in $\langle x,a_{1}%
,a_{2},a_{3},a_{4}\rangle$ are zero provided $p\neq2,3,5$.

\item All Lie products of multidegree $(5,1,1,1,1,1)$ in $\langle
x,a_{1},a_{2},a_{3},a_{4},a_{5}\rangle$ are zero provided $p\neq2,3,5,7,31$.

\item All Lie products of multidegree $(4,1,1,1,1,1,1)$ in $\langle
x,a_{1},a_{2},a_{3},a_{4},a_{5},a_{6}\rangle$ are zero provided $p\neq2,3,5 $.

\item All Lie products of multidegree $(3,1,1,1,1,1,1)$ in $\langle
x,a_{1},a_{2},a_{3},a_{4},a_{5},a_{6}\rangle$ are zero provided $p\neq2,3,5,7$.

\item All Lie products of multidegree $(2,1,1,1,1,1,1)$ in $\langle
x,a_{1},a_{2},a_{3},a_{4},a_{5},a_{6}\rangle$ are zero provided $p\neq2,3$.
\end{itemize}

Clearly any Lie product of multidegree $(6,1)$ in the subalgebra $\langle
x,a_{1}\rangle$ of $L$ is zero (for all primes $p$). Also any Lie product of
multidegree $(1,1,1,1,1,1,1)$ in $\langle x,a_{1},a_{2},a_{3},a_{4}%
,a_{5},a_{6}\rangle$ is zero for all primes $p>5$. So these results establish
the last line of the table above.

The computations needed to prove these seven claims are all similar, and so I
will just describe the computation needed to show that all Lie products of
multidegree $(4,1,1,1,1,1,1)$ in $\langle x,a_{1},a_{2},a_{3},a_{4},a_5,a_6\rangle$ 
are zero provided $p\neq2,3,5$.

So let $M$ be the Lie algebra over $\mathbb{Q}$ generated by $x,a_{1}%
,a_{2},\ldots ,a_{6}$ satisfying the relations 
$[a_{i},a_{j}]=0$ for all $i,j$,
and also satisfying the condition that any Lie product of the generators in
$M$ with multidegree $(d_{x},d_{1},d_{2},\ldots ,d_{6})$ is zero if $d_{x}>4$ or
$d_{i}>1$ for some $i$. Note that we do {\emph not} enforce the 5-Engel
identity in $M$. Then my Magma program shows that $M$ has class 10 and
dimension 5705. The program gives a basis for $M$ consisting of left-normed Lie
products of the generators $x,a_{1},a_{2},\ldots ,a_{6}$, and gives the product
of any two basis elements as a linear combination of the basis elements.
Fortunately the structure constants arising as coefficients in these linear
combinations are all integers. (Perhaps this is just lucky, or perhaps not!
Who knows? The structure constants arising in all these seven computations are all
integers.) Since the structure constants are all integers, the presentation
for $M$ can be used to define a Lie algebra $M_{K} $ over any field $K$, and
this Lie algebra has dimension 5705 over $K$ and satisfies all the relations
of $M$. But we have to consider the possibility that for some fields $K$ of
finite characteristic the largest Lie algebra over $K$ on
these generators, satisfying these relations, might have dimension bigger that
5705. We can discount this possibility as follows. Let $A$ be any Lie algebra
generated by elements $x,a_{1},a_{2},\ldots ,a_{6}$ satisfying the relations of
$M$. Then $[A,A]$ is generated by left-normed Lie products of the form%
\[
\lbrack x,a_{i_{1}},a_{i_{2}},\ldots,a_{i_{r}},x,\ldots,x]
\]
which have degree at most 4 in $x$, and where $1\leq r\leq 6$ and $1\leq
i_{1}<i_{2}<\ldots<i_{r}\leq 6$. It is straightforward to count the number of
basic Lie products in these generators of $[A,A]$ which have multidegree
$(d_{x},d_{1},d_{2},\ldots ,d_{6})$ with $d_{x}\leq4$ and with $d_{i}=0$ or 1
for $i=1,2,\ldots 6$. This enables us to obtain an upper bound on the dimension of
$A$, and this upper bound is 5705.

An alternative way of showing that $M_K$  is the largest Lie algebra over $K$
on generators $x,a_{1},a_{2},\ldots ,a_{6}$ which satisfies the relations of $M$
is as follows. Let $K$ be any field, and let $A$ be the abelian Lie algebra
over $K$ with basis $a_1,a_2,\ldots ,a_6$. Let ${\cal N}_4$ be the variety of Lie
algebras over $K$ which are nilpotent of class 4. Then define $B$ to be the 
wreath ${\cal N}_4$-product of $\langle x \rangle $ and $A$. (See [2, Section 3.7].)
Let $X$ be the ideal of $B$ generated by $x$. Then $B=A\oplus X$ and $X$ is 
a free nilpotent of class 4 Lie algebra over $K$ with free generators
\[
[x,a_i,a_j,\ldots ,a_k]
\]
with $1 \le i \le j \le \ldots \le k \le 6$. We then let $I$ be the ideal of $B$
spanned by all Lie products in the generators $x,a_1,a_2,\ldots ,a_6$ which have
multidgree $(d_x,d_1,d_2,\ldots ,d_6)$ with $d_i>1$ for some $i \in \{1,2,\ldots ,6\}$.
Then $B/I$ is the largest Lie algebra over $K$ generated by 7 elements satisfying
the relations of $M$, and the dimension of $B/I$ is the same for all fields $K$.

So we can assume that $M_{K}$ is the largest Lie algebra over $K$ on
generators $x,a_{1},a_{2},\ldots ,a_{6}$ which satisfies the relations of $M$.
Next, we generate all 5-Engel relations of weight 10 in the generators of $M$.
Specifically, we evaluate
\[
\sum_{\sigma\in\text{Sym}(5)}[y,x_{1\sigma},x_{2\sigma},x_{3\sigma}%
,x_{4\sigma},x_{5\sigma}]
\]
for all possible sequences $y,x_1,x_2,x_3,x_4,x_5$ of basis elements of $M$
whose multidegrees sum to $(4,1,1,1,1,1,1)$.
The number of basis elements of $M$ of weight 10 is $1024$, and we suppose
that these are $b_{1},b_{2},\ldots,b_{1024}$. All these basis elements of weight
10 have multidegree $(4,1,1,1,1,1,1)$. The 5-Engel relations are
all of the form $\sum_{i=1}^{1024}m_{i}b_{i}=0$ where the coefficients $m_{i}$
are all integers. We collect all the vectors 
\[
(m_{1},m_{2},\ldots,m_{1024})
\]
from these relations in an integer matrix with $1024$ columns. Then we check
that the rank of the matrix is $1024$, and that the elementary divisors of the
matrix only involve primes $2,3,5$. This establishes our claim that all Lie
products in $L$ of multidegree $(4,1,1,1,1,1,1)$ in $\langle x,a_{1},a_{2}%
,\ldots ,a_{6}\rangle$ are zero provided $p\neq2,3,5$.

\end{document}